\documentclass[12pt]{article}
\usepackage{amsmath,amssymb,amsfonts}
\usepackage{color}
\usepackage{pstcol,pstricks}
\usepackage{graphicx}

\numberwithin{equation}{section}

{ \pagestyle{plain}

\def\cwedge{\bigcirc\kern-1.07em\wedge\ }

\newcommand{\qed}{\hfill\fbox{}\par\vspace{.2cm}}

\begin{document}

%
%

\begin{center}
{\LARGE {{A one-parameter family of totally umbilical hyperspheres
in the nearly {K\"ahler} $6$-sphere}}}
\end{center}
{
\begin{center}
{\large Jihong Bae$^{\dag}$, JeongHyeong Park$^{\dag}$, Kouei
Sekigawa$^{\ddag}$}
\end{center} }
\begin{center}
$^{\dag}$Sungkyunkwan University,
 Suwon, Korea\\
$^{\ddag}$Niigata University,
    Niigata, Japan\end{center}

\begin{abstract}
{{We discuss two kinds of almost contact metric structures on a
one-parameter family of totally umbilical hyperspheres in the nearly
{K\"ahler} unit $6$-sphere $S^{6}$.}}
\end{abstract}
\noindent {\it Mathematics Subsect Classification (2010)} : 53C25, 53D10

\section{Introduction}\label{sec1}
An odd dimensional smooth manifold $M$ with a quadruple $(\phi, \xi, \eta, g)$ of a $(1,1)$-tensor field $\phi$, a vector field $\xi$, a 1-form $\eta$ and a Riemannian metric $g$ satisfying the following conditions is called an almost contact metric manifold;
\begin{equation}\label{11}
\begin{split}
\phi^2=&-I+\eta\otimes\xi, \quad \eta(\xi)=1\\ &\phi\xi =0,\quad
\eta\circ\phi=0
\end{split}
\end{equation}
and
\begin{equation}\label{12}
\begin{split}
g(\phi{X},\phi{Y})&=g(X,Y)-\eta(X)\eta(Y),
\end{split}
\end{equation}
for any $X, Y \in\mathfrak{X}(M)$, where $\mathfrak{X}(M)$ denotes the
Lie algebra of all smooth vector fields on $M$. Further, an almost contact metric manifold $(M, \phi, \xi, \eta, g)$ is called a contact metric manifold if it satisfies the following condition;
\begin{equation}\label{13}
\begin{split}
d\eta(X,Y)=g(X,\phi{Y}),
\end{split}
\end{equation}
for any $X, Y \in\mathfrak{X}(M)$. In \cite{KPS}, authors defined a
new class of almost contact metric manifolds, say, the class of quasi
contact metric manifolds which are a generalization of contact metric manifolds, and the basic properties for quasi contact metric manifolds also have been obtained \cite{CKPSS, PSS}. Further, the authors raised the following question based on the discussion.\\
\\
{\bf{Question A}.} Does there exist a $(2n+1)(\geq 5)$-dimensional quasi contact metric
manifold which is not a contact metric manifold?\\
\\
Concerning the above Question A, authors discussed oriented hypersurfaces in a quasi {K\"ahler} manifold which are quasi contact metric manifolds
with respect to the naturally induced almost contact metric structure, and obtained the following results in \cite{BPS}.\\

\noindent {\bf{Theorem B}.} Let $\bar{M}=(\bar{M},\bar{J},\bar{g})$
be a nearly K\"{a}hler manifold and $M$ be a
hypersurface of $\bar{M}$ oriented by a unit normal
vector field $\nu$. Then $M = (M,\phi,\xi,\eta, g)$ is
a quasi contact metric manifold with respect to the
naturally induced almost contact metric structure
$(\phi,\xi,\eta,g)$ if and only if it satisfies the equality
\begin{equation*}
\begin{split}
g((A\phi +\phi A)X, Y) = -2g(\phi X, Y)
\end{split}
\end{equation*}
for any $X,Y\in \mathfrak{X}(M)$, where $A$ is the
shape operator with respect to the unit normal vector
field $\nu$, and hence, $M$ is a contact metric manifold.
\\

\noindent{\bf{Theorem C}.} There does not exist oriented totally umbilical hypersurface in the nearly {K\"ahler} unit $6$-sphere which is a quasi contact metric manifold with respect to the naturally induced almost contact metric structure.\\
\\
In the present paper, we provide explicit examples of totally
umbilical hypersurfaces in the nearly {K\"ahler} unit $6$-sphere
which support Theorem B and Theorem C.

\section{Preliminaries}\label{sec2}
First, we shall recall fundamental the nearly {K\"ahler} structure
on a unit $6$-sphere $S^{6}$. Let $\mathfrak{C}$ be the Cayley
algebra $\mathfrak{C}=\{x=x_{0}+\sum_{i=1}^{7}x_{i}e_{i}~|~x_{0},
x_{i}\in \mathbb{R}, e_{i}^{2}=-1~(1\leq i \leq7)\}$, and
$\mathfrak{C}_{+}=\{x = \sum_{i=1}^{7}x_{i}e_{i} \in
\mathfrak{C}~|~x_{i}\in \mathbb{R}~(1\leq i \leq 7)\}$ both set of
all pure imaginary Cayley numbers. Here, the multiplication
operation on $\mathfrak{C}$ is defined by the figure below;
\begin{figure}\label{21}
\centering
\includegraphics[width=8cm]{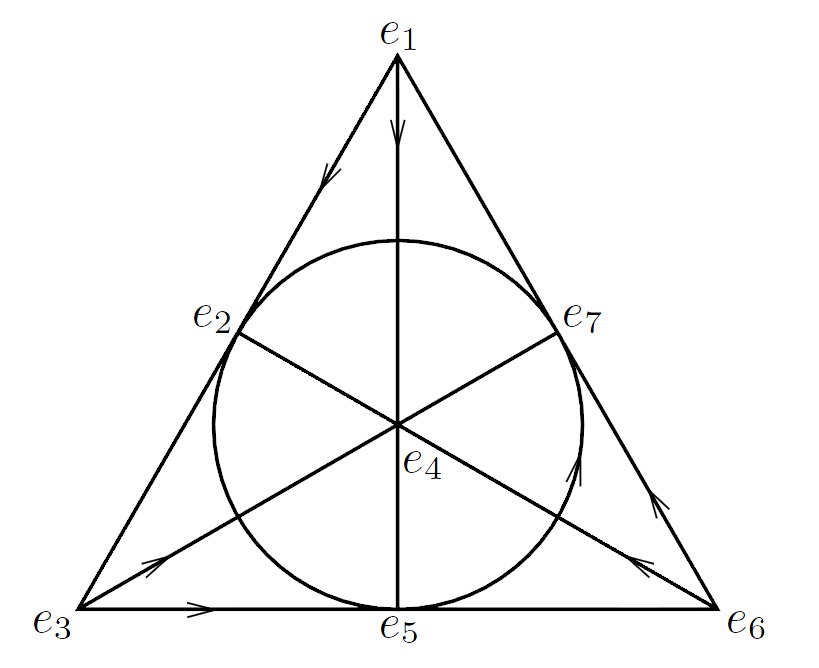}\\
\caption{} \label{21}
\end{figure}

 We denote by $<,>$ the canonical inner product on
$\mathfrak{C}$ and let $|x|=\sqrt{<x,x>}$ (the length of
$x\in\mathfrak{C}$). Then, $(\mathfrak{C}, <,>)$
(resp. $(\mathfrak{C}_{+}, <,>)$) can be identified with $8$-dimensional Euclidean space $\mathbb{E}^{8}$ (resp. $7$-dimensional Euclidean space $\mathbb{E}^{7}$) in the natural
way. We also define cross product $x \times y$ for $x, y
\in\mathfrak{C}_{+}$ by { $x \times y=xy+<x, y>1 (\in\mathfrak{C}_{+})$}. Here, we identify
$e_{i}\in\mathfrak{C}_{+}~(1 \leq i \leq 7)$ with the coordinate
vector field $\frac{\partial}{\partial x_{i}}$ (denoted by
$\partial_{i}$ briefly) in our arguments and adopt them
alternatively in the forthcoming arguments. We denote by D the
Levi-Civita connection on $\mathbb{E}^{7}$ with respect to the
Riemannian metric induced from the inner product $<,>$. Let $S^{6}$
be a unit $6$-sphere in
$\mathbb{E}^{7}(\simeq\mathfrak{C}_{+})$ centered at the origin $o$.
Then, $S^{6}$ is expressed as
$S^{6}=\{x\in\mathfrak{C}_{+}~|~|x|=1\}$.

For any point $x\in S^{6}$, we denote by $N_{x}$ the outward
oriented unit normal vector with initial point $x$,
$N_{x}=\overrightarrow{ox}$. In this paper, we identify $N_{x}(x\in
S^{6})$ with the position vector $x(\in\mathfrak{C}_{+})$. The unit
normal vector $N$ is also written as
$N=\sum_{i=1}^{7}x_{i}\partial_{i}$ in terms of the coordinate vector
fields $\partial_{i}(1\leq i \leq 7)$. Here we note that the tangent
space $T_{x}S^{6}$ can be regarded as the subspace
$\{y\in\mathfrak{C}_{+}~|~<y,x>=0\}$ of $\mathfrak{C}_{+}$. Now, we
define $(1,1)$-tensor field $J$ on $S^{6}$ by
\begin{equation}\label{22}
\begin{split}
J_{x}y=N_{x}\times y=x\times y(=xy),\quad y\in T_{x}S^{6}.
\end{split}
\end{equation}
Then, we may easily check that $J$ is an almost complex structure on
$S^{6}$ and $(J,\bar{g})$ is a nearly K\"ahler structure on $S^{6}$,
namely, $(\bar{\nabla}_{X}J)Y=-(\bar{\nabla}_{Y}J)X$ holds for any
vector fields $X,Y$ tangent to $S^{6}$, where $\bar{g}$ and
$\bar{\nabla}$ are the Riemannian metric on $S^{6}$ induced from the
inner product $<,>$ on $\mathfrak{C}_{+}$ and $\bar{\nabla}$ is the
Levi-Civita connection of $\bar{g}$, respectively. We shall call the
nearly {K\"ahler} structure $(J, \bar{g})$ given above on $S^{6}$
the standard one.

\section{One parameter family of totally umbilical hyperspheres in $S^{6}$}\label{sec3}
First, for each real number $r~(-1<r<1)$, we define hypersurface $M_{r}$ by
\begin{equation*}
\begin{split}
M_{r}&=S^{6}\cap\{x=\sum_{i=1}^{6}x_{i}e_{i}+re_{7}\in\mathfrak{C}_{+}|~x_{i}\in\mathbb{R}~(1\leq i\leq 6)\}\\
&=\{x=\sum_{i=1}^{6}x_{i}e_{i}+re_{7}\in\mathfrak{C}_{+}|~\sum_{i=1}^6x_{i}^{2}=1-r^{2}\}.
\end{split}
\end{equation*}
We observe that $M_{r}$ is diffeomorphic to a $5$-sphere $S^{5}$.

Now, let $x$ be any point of $M_{r}$ and $\gamma_{x}$ be the smooth curve in $M_{r}$ through $x=\gamma_{x}(\theta)~(0<\theta<\pi)$ defined by
\begin{equation}\label{31}
\begin{split}
\gamma_{x}(t)=(\cos t)e_{7}+(\frac{1}{\sqrt{1-r^{2}}}\sin t)\sum_{i=1}^{6}x_{i}e_{i}~(0\leq t \leq \pi),
\end{split}
\end{equation}
where $\cos\theta=r,~\sin\theta=\sqrt{1-r^{2}}$. We here define a vector field $\nu$ on $M_r$ by
\begin{equation}\label{32}
\begin{split}
\nu_{x}&=\frac{d}{dt}|_{t=\theta}\gamma_{x}(t)\\&=-(\sin\theta)e_{7}+(\frac{1}{\sqrt{1-r^{2}}}\cos\theta)\sum_{i=1}^{6}x_{i}e_{i}
\\&=\frac{r}{\sqrt{1-r^{2}}}\sum_{i=1}^{6}x_{i}e_{i}-\sqrt{1-r^{2}}e_{7}.
\end{split}
\end{equation}
Thus, from \eqref{32}, the following equalities hold  for any $x\in M_{r}$;\\
\begin{equation}\label{33}
\begin{split}
\bar{g}(\nu_{x},\nu_{x})&=<\nu_{x},\nu_{x}>=1,\\
<\nu_{x},N_{x}>=<\nu_{x},x>=&\frac{r}{\sqrt{1-r^{2}}}(1-r^{2})-r\sqrt{1-r^{2}}=0.
\end{split}
\end{equation}
On the other hand, for any $x=\sum_{i=1}^{6}x_{i}e_{i}+re_{7}\in M_{r}$, we may find an integer
$a$ $(1\leq a\leq 6)$ such that $x_{a}\neq 0$ and fix it. Now, we
shall define a smooth curve $\alpha_{a,b}(s)$ $(1\leq b \leq 6,~ b\ne a)(-\pi<s<\pi) $ through the point $x = \alpha_{a,b}(0)$ by
\begin{equation}\label{34}
\begin{split}
\alpha_{a,b}(s)=(\sqrt{x_{a}^{2}+x_{b}^{2}}\cos(s+\theta_{a,b}))e_{a} +(\sqrt{x_{a}^{2}+x_{b}^{2}}\sin(s+\theta_{a,b}))e_{b}+{
{\sum_{1 \leq i \leq 6,i\neq a,b}x_ie_i+re_7}},
\end{split}
\end{equation}
where $\cos\theta_{a,b}=\frac{x_{a}}{\sqrt{x_{a}^{2}+x_{b}^{2}}},~\sin\theta_{a,b}=\frac{x_{b}}{\sqrt{x_{a}^{2}+x_{b}^{2}}}~(0\leq\theta_{a,b}<2\pi)$.
Then, from \eqref{34}, we have
\begin{equation}\label{35}
\begin{split}
\frac{d}{ds}|_{s=0}\alpha_{a,b}(s)&=-(\sqrt{x_{a}^{2}+x_{b}^{2}}\sin\theta_{a,b})e_{a}+(\sqrt{x_{a}^{2}+x_{b}^{2}}\cos\theta_{a,b})e_{b}\\&=-x_{b} e_{a} + x_{a} e_{b}~(=- x_{b} \partial_{a}+x_{a}\partial_{b})
\end{split}
\end{equation}
at $x$. We here set
\begin{equation}\label{36}
\begin{split}
X_{a,b}=-x_{b}e_{a}+x_{a}e_{b}~(=-x_{b}\partial_{a}+x_{a}\partial_{b}).
\end{split}
\end{equation}
From \eqref{35} and \eqref{36}, it follows that
\begin{equation}\label{37}
\begin{split}
T_{x}M_{r}=span_{\mathbb{R}}\{X_{a,b}~(b\neq a,~1\leq b \leq 6)\}
\end{split}
\end{equation}
and
\begin{equation}\label{38}
\begin{split}
\bar{g}(X_{a,b},\nu_{x})=0,
\end{split}
\end{equation}
at $x \in M_{r}$. Thus, from \eqref{37} and \eqref{38}, we can see that $\nu_{x}$ is a unit normal vector at any $x \in M_{r}$ in $S^{6}$, namely the vector field $\nu$ is a unit normal vector field on $M_{r}$ in $S^{6}$.

Now, since $S^{6}$ is a totally umbilical hypersurface in
$\mathbb{E}^{7}~(\backsimeq \mathbb{C}_{+})$ with respect to the
unit normal vector field $N$, the corresponding shape operator $\bar{A}$ is
given by $\bar{A}=-I$. Thus, taking account of the Gauss formula, we
have
\begin{equation}\label{39}
\begin{split}
D_{X_{a,b}}\nu=\bar{\nabla}_{X_{a,b}}\nu.
\end{split}
\end{equation}
From \eqref{32}, the unit normal vector field $\nu$ can be expressed by
\begin{equation}\label{310}
\begin{split}
\nu=\frac{r}{\sqrt{1-r^{2}}}\sum_{i=1}^{6} x_{i}\partial_{i} -\sqrt{1-r^{2}}\partial_{7}.
\end{split}
\end{equation}
Thus, from \eqref{36}, \eqref{39} and \eqref{310}, we have
\begin{equation}\label{311}
\begin{split}
D_{X_{a,b}}\nu&=\frac{r}{\sqrt{1-r^{2}}}(-x_{b}\partial_{a}+x_{a}\partial_{b})=\frac{r}{\sqrt{1-r^{2}}}X_{a,b}
\end{split}
\end{equation}
for any $X_{a,b}$ at any point $x \in M_{r}$. Therefore, from \eqref{311}, we see that $(M_{r},g)$ is a totally umbilical hypersurface of $(S^{6},\bar{g})$ with the shape operator $A=-\frac{r}{\sqrt{1-r^{2}}}I$ with respect to the unit normal vector field $\nu$ on $(M_{r},g)$ in $(S^{6},\bar{g})$.
\section{Almost contact metric structures on $(M_r, g)$}\label{sec4}
In this section, we define two kinds almost contact metric
structures on $(M_{r}, g)$ and discuss their respective geometric
properties. \noindent First, let $\xi$ be the unit vector field on
$M_{r}$ defined by
\begin{equation}\label{41}
\begin{split}
\xi = -J\nu = -N\times \nu.
\end{split}
\end{equation}
Then, from \eqref{41}, it follows that the vector field
$\xi$ is orthogonal to both of the vector fields $N$ and $\nu$ along
$M_{r}$. Further, from Fig.\eqref{21}, \eqref{32} and \eqref{41}, we have
\begin{equation}\label{42}
\begin{split}
\xi &= -(\sum_{i=1}^6 x_ie_i+re_7)\times
(\frac{r}{\sqrt{1-r^2}}\sum_{j=1}^6 x_je_j-\sqrt{1-r^2}e_7)\\
&=\sqrt{1-r^2}(\sum_{i=1}^6 x_ie_i)\times e_7-\frac{r^2}{\sqrt{1-r^2}}e_7
\times (\sum_{j=1}^6 x_je_j)\\
&=\frac{1}{\sqrt{1-r^2}}(x_6e_1+x_5e_2+x_4e_3-x_3e_4-x_2e_5-x_1e_6).
\end{split}
\end{equation}
From \eqref{42}, $\xi$ is also rewritten as
\begin{equation}\label{43}
\begin{split}
\xi
=\frac{1}{\sqrt{1-r^2}}(x_6\partial_1+x_5\partial_2+x_4\partial_3-x_3\partial_4-x_2\partial_5-x_1\partial_6).
\end{split}
\end{equation}
Thus, the 1-form $\eta$ dual to the vector field $\xi$ is given by
\begin{equation}\label{44}
\begin{split}
\eta
=\frac{1}{\sqrt{1-r^2}}(x_6dx_1+x_5dx_2+x_4dx_3-x_3dx_4-x_2dx_5-x_1dx_6).
\end{split}
\end{equation}
From \eqref{44}, we also have
\begin{equation}\label{45}
\begin{split}
d\eta =-\frac{2}{\sqrt{1-r^2}}(dx_1\wedge dx_6+dx_2\wedge
dx_5+dx_3\wedge dx_4).
\end{split}
\end{equation}
From \eqref{44} and \eqref{45}, we have further
\begin{equation}\label{46}
\begin{split}
\eta \wedge (d\eta)^2 =-\frac{8}{(\sqrt{1-r^2})^3}
\{-&x_1dx_2\wedge dx_3\wedge dx_4\wedge dx_5\wedge dx_6\\
+&x_2dx_1\wedge dx_3\wedge dx_4\wedge dx_5\wedge dx_6\\
-&x_3dx_1\wedge dx_2\wedge dx_4\wedge dx_5\wedge dx_6\\
+&x_4dx_1\wedge dx_2\wedge dx_3\wedge dx_5\wedge dx_6\\
-&x_5dx_1\wedge dx_2\wedge dx_3\wedge dx_4\wedge dx_6\\
+&x_6dx_1\wedge dx_2\wedge dx_3\wedge dx_4\wedge dx_5 \}\neq 0.
\end{split}
\end{equation}
Therefore, $\eta$ is a contact form on $M_{r}$. Now, we shall
show $\nabla_{\xi}\xi=0$. From \eqref{32}, we have
\begin{equation}\label{47}
\begin{split}
D_{\xi}\xi=-\frac{1}{1-r^2}\sum_{i=1}^{6}x_i\partial_i
\end{split}
\end{equation}
Taking account of the Gauss formula for $(S^6, \bar{g})$ and $(\mathbb{E}^7,<,
>)$, we have
\begin{equation}\label{48}
\begin{split}
D_{\xi}\xi&=\bar{\nabla}_{\xi}\xi-N =\bar{\nabla}_{\xi}\xi-\sum_{i=1}^6 x_i\partial_i-r\partial_7.
\end{split}
\end{equation}
On the other hand, since $(M_{r}, g)$ is a totally umbilical
hypersurface of $(S^6, \bar{g})$ with the shape operator
$A=-\frac{r}{\sqrt{1-r^2}}{I}$ with respect to the unit normal vector field
$\nu$, from \eqref{310}, taking account of the Gauss formula, we get
\begin{equation}\label{49}
\begin{split}
\bar{\nabla}_{\xi}\xi
&={\nabla}_{\xi}\xi-\frac{r}{\sqrt{1-r^2}}\nu\\
&={\nabla}_{\xi}\xi-\frac{r}{\sqrt{1-r^2}}(\frac{r}{\sqrt{1-r^2}}\sum_{i=1}^6
x_i\partial_i-\sqrt{1-r^2}\partial_7)\\
&={\nabla}_{\xi}\xi-{\frac{r^2}{1-r^2}}\sum_{i=1}^6
x_i\partial_i+r\partial_7.
\end{split}
\end{equation}
Then, from \eqref{47}$\thicksim$\eqref{49}, we have
\begin{equation*}
\begin{split}
-\frac{1}{1-r^2}\sum_{i=1}^6
x_i\partial_i=\nabla_{\xi}\xi-(1+\frac{r^2}{1-r^2})\sum_{i=1}^6
x_i\partial_i,
\end{split}
\end{equation*}
and hence
\begin{equation}\label{410}
\begin{split}
{\nabla}_{\xi}\xi=0.
\end{split}
\end{equation}
From \eqref{410}, it follows that each integral curve of the vector
field $\xi$ is a geodesic of $(M_{r}, g)$. Thus, taking account of the
definition of the vector field $\xi$ in \eqref{41}, we see that
$(M_{r}, g, \xi)$ is a Hopf hypersurface in $(S^6, J, \bar{g})$.
Further, since $(M_{r}, g)$ is a totally umbilical hypersurface in
$(S^6, \bar{g})$ with the shape operator
$A=-\frac{r}{\sqrt{1-r^2}}I$, from the Gauss equation for $(M_{r},
g)$, we see that the curvature tensor $R$ of $(M_{r}, g)$ is given
\begin{equation}\label{411}
\begin{split}
R(X, Y)Z&=g(Y, Z)X-g(X, Z)Y+\frac{r^2}{1-r^2}(g(Y, Z)X-g(X, Z)Y)\\
&=\frac{1}{1-r^2}(g(Y, Z)X-g(X, Z)Y),
\end{split}
\end{equation}
for any $X, Y, Z \in {T_{x}M_{r}}$. From \eqref{411}, it follows that
$(M_{r}, g)$ is a hypersurface of $(S^6, \bar{g})$ of constant
sectional curvature $\frac{1}{1-r^2}$. We define $(1,1)$-tensor
field $\phi$ on $M_{r}$ by
\begin{equation}\label{412}
\begin{split}
\phi X=JX-\eta (X)\nu
\end{split}
\end{equation}
for any $X \in {T_{x}M_{r}}$. Then, from \eqref{41},
\eqref{44} and \eqref{412}, we see that $(\phi, \xi, \eta, g)$ is
the naturally induced almost contact metric structure on $M_{r}$. Now,
choose $x=\sum_{i=1}^6 x_ie_i+re_7 \in M_{r}$ arbitrary.
Without loss of essentiality, we may suppose $x_1\neq 0$, for example.
Then, from \eqref{22}, \eqref{36} and \eqref{412}, taking account of Fig.\eqref{21}, we
have
\begin{equation}\label{413}
\begin{split}
X_{1, 2}=-x_2\partial_1+x_1\partial_2, \quad X_{1,
3}=-x_3\partial_1+x_1\partial_3,
\end{split}
\end{equation}
and
\begin{equation}\label{414}
\begin{split}
\phi X_{1, 3}&=JX_{1, 3}-\eta(X_{1, 3})\nu\\
&=(x_1x_2-\frac{r}{1-r^2}(x_1^2x_4-x_1x_3x_6))\partial_1\\
&+(-(x_1^2+x_3^2)-\frac{r}{1-r^2}(x_1x_2x_4-x_2x_3x_6))\partial_2\\
&+(x_2x_3-\frac{r}{1-r^2}(x_1x_3x_4-x_3^2x_6))\partial_3\\
&+((-x_3x_5+rx_1)-\frac{r}{1-r^2}(x_1x_4^2-x_3x_4x_6))\partial_4\\
&+(x_1x_6+x_3x_4-\frac{r}{1-r^2}(x_1x_4x_5-x_3x_5x_6))\partial_5\\
&+(-(rx_3+x_1x_5)-\frac{r}{1-r^2}(x_1x_4x_6-x_3x_6^2))\partial_6.
\end{split}
\end{equation}
Thus, from \eqref{413} and \eqref{414}, we have
\begin{equation}\label{415}
\begin{split}
g(X_{1, 2}, \phi X_{1, 3})=-x_1(x_1^2+x_2^2+x_3^2) (\neq 0).
\end{split}
\end{equation}
On the other hand, from \eqref{45} and \eqref{413}, we have
\begin{equation}\label{416}
\begin{split}
d\eta(X_{1, 2}, X_{1, 3})&=-\frac{2}{\sqrt{1-r^2}}(dx_1\wedge
dx_6+dx_2\wedge dx_5+dx_3\wedge
dx_4)(-x_2\partial_1+x_1\partial_2, -x_3\partial_1+x_1\partial_3)\\
&=0.
\end{split}
\end{equation}
Thus, from \eqref{415} and \eqref{416}, we have
\begin{equation}\label{417}
\begin{split}
d\eta(X_{1, 2}, X_{1, 3})\neq g(X_{1, 2}, \phi X_{1, 3}).
\end{split}
\end{equation}
We may also derive the similar conclusion as \eqref{417} for the
other cases $x_b \neq 0\\ (3\leq b \leq 6)$. Thus, from
\eqref{417}, the almost contact metric manifold $(M_{r},\phi, \xi, \eta, g)$ is not a contact metric manifold for any $r~(-1<r<1)$. This supports Theorem C.\\
Now, let $\tau$ be the scalar curvature of $(M_{r}, g)$ and define the
smooth functions $f$ on $M_r$ and the mean curvature $\alpha$ respectively by
\begin{equation}\label{418}
\begin{split}
f=g(A\xi, \xi)
\end{split}
\end{equation}
and
\begin{equation}\label{419}
\begin{split}
\alpha=\frac{1}{5}trA.
\end{split}
\end{equation}
Then, from \eqref{411}, since $A=-\frac{r}{\sqrt{1-r^2}}I$ and
$g(\xi, \xi)=1$, we have
\begin{equation}\label{420}
\begin{split}
\tau=\frac{20}{1-r^2}, \quad f=-\frac{r}{\sqrt{1-r^2}}, \quad
\alpha=-\frac{r}{\sqrt{1-r^2}}.
\end{split}
\end{equation}
Then, from \eqref{420}, we may check that the hypersurface $(M_{r}, g, \xi,\eta,g)$ satisfies
the following equality;
\begin{equation}\label{421}
\begin{split}
\tau = 20+5\alpha(5\alpha-f),
\end{split}
\end{equation}
for any real number $r~(-1<r<1)$.
Here, we note that $(M_{0},\phi, \xi, \eta, g)$ is totally geodesic in $(S^6, g)$
if and only if $r=0.$ This shows that
the statement of the result(\cite{DA}, Theorem 1.1) is inadequate
by taking account of the 1-parameter family of the Hopf
hypersurfaces $(M_r, g, \xi)$ $(-1<r<1)$ of the nearly {K\"ahler}
6-sphere $(S^6, J, \bar{g})$.

\noindent
{{Now,}} we define another almost contact metric structure on the
hypersurface $(M_{r}, g)$ and discuss on the geometric properties.
Let $\phi'$ be the $(1,1)$-tensor field on $M_{r}$ defined by
\begin{equation}\label{422}
\begin{split}
\phi '\xi=0
\end{split}
\end{equation}
and
\begin{equation}\label{423}
\begin{split}
\phi 'X=-\nu \times X=-\nu X
\end{split}
\end{equation}
for any $X\in\mathfrak{X}(M_{r})$ with $X \perp \xi$. Then, from {\eqref{41}}, \eqref{422} and \eqref{423}, we may check
\begin{equation}\label{424}
\begin{split}
<\phi 'X, N>&=-<\nu X, N>= -<\nu, \nu><\nu X, N>\\
&=-<\nu(\nu X), \nu N>=-<\nu^{2} X, \nu N >\\
&=<X, \nu N >= -<X, N\nu >\\
&=<X, \xi>=0,
\end{split}
\end{equation}
\begin{equation}\label{425}
\begin{split}
<\phi 'X, \nu>&=<-\nu X, \nu>\\
&=-<\nu, \nu><\nu X, \nu>\\
&=-<\nu^{2} X,\nu^{2}>=<X,1>\\
&=0.
\end{split}
\end{equation}
\begin{equation}\label{426}
\begin{split}
<\phi 'X, \xi>&=<-\nu X, -N\nu>\\
&=<-\nu X, \nu N>\\
&=-<\nu, \nu><X, N>\\
&=0.
\end{split}
\end{equation}
Thus, from \eqref{423} $\sim$ \eqref{426}, we have finally
\begin{equation}\label{427}
\begin{split}
\phi '^2X =\phi '(\phi ' X)=\nu(\nu X)=\nu^{2}X =-X.
\end{split}
\end{equation}
Taking account of \eqref{422} $\sim$ \eqref{427}, we see that $(\phi', \xi, \eta, g)$ is an almost contact metric structure on $M_{r}$.
Further, we may note that the almost contact metric manifold $(M_0,\phi ', \xi, \eta, g)$ coincides with the almost contact metric
manifold introduced in (\cite{B1}, p. 64) which is different from the naturally induced cone from the nearly K\"ahler structure $(J,\bar{g})$ on $S^6$ with respect to the unit normal vector field
$\nu$. Later, we shall show that $(\phi ', \xi, \eta, g)$ is a contact metric structure on the hypersurface $M_0$.

\noindent {For any $X$ $\in\mathfrak{X}(M_{r})$, we set}
\begin{equation}\label{428}
\begin{split}
Y=X-\eta(X)\xi.
\end{split}
\end{equation}
Then, $Y$ $\in\mathfrak{X}(M_{r})$ and $Y \perp \xi$. From
\eqref{41}, \eqref{422} and \eqref{423}, we have
\begin{equation}\label{429}
\begin{split}
\phi 'Y&=-\nu \times Y =-\nu \times (X-\eta(X)\xi)\\
&=-\nu \times X + \eta(X)(\nu \times \xi)=-\nu X + \eta(X)\nu\xi\\
&=-\nu X-\eta(X)\nu(N\nu) = -\nu X + \eta(X)\nu^{2}N \\
&=-\nu X - \eta(X)N,
\end{split}
\end{equation}
for any $X$ $\in\mathfrak{X}(M_{r})$. Comparing \eqref{412} and
\eqref{429}, we see that $\phi X \neq \phi 'X$ for $X$
$\in\mathfrak{X}(M_{r})$ with $X \perp \xi$. It is known that the
almost contact metric structure $(\phi ', \xi, \eta, g)$ on the
hypersphere $M_0$ in the nearly {K\"ahler} 6-sphere $S^6$ is a
contact metric structure by (\cite{B1}, p. 64). Here, we shall
provide an exact proof for this fact, now, we choose a point
$x=\sum_{i=1}^6 x_ie_i \in M_0$ arbitrary and fix it. Here, for
our purpose without also discuss in the case where $x_i \neq 0$,
now, we set
\begin{equation}\label{430}
\begin{split}
Y_{1, b}=X_{1, b}-\eta(X_{1, b})\xi \quad(1<b\leq 6)
\end{split}
\end{equation}
for any $X$ $\in\mathfrak{X}(M_{r})$. Then from \eqref{430}, taking
account of \eqref{36}, \eqref{42} with $r=0$, \eqref{43} and Fig.\eqref{21},
we have
\begin{equation}\label{431}
\begin{split}
Y_{1,2}&=(-x_2+x_2x_6^2-x_1x_5x_6)\partial_1+(x_1+x_2x_5x_6-x_1x_5^2)\partial_2\\
&+(x_2x_4x_6-x_1x_4x_5)\partial_3+(-x_2x_3x_6+x_1x_3x_5)\partial_4\\
&+(-x_2^2x_6+x_1x_2x_5)\partial_5+(-x_1x_2x_6+x_1^2x_5)\partial_6,
\end{split}
\end{equation}
\begin{equation*}
\begin{split}
Y_{1,3}&=(-x_3+x_3x_6^2-x_1x_4x_6)\partial_1+(x_3x_5x_6-x_1x_4x_5)\partial_2\\
&+(x_1+x_3x_4x_6-x_1x_4^2)\partial_3+(-x_3^2x_6+x_1x_3x_4)\partial_4\\
&+(-x_2x_3x_6+x_1x_2x_4)\partial_5+(-x_1x_3x_6+x_1^2x_4)\partial_6,
\end{split}
\end{equation*}
\begin{equation*}
\begin{split}
Y_{1,4}&=(-x_4+x_4x_6^2+x_1x_3x_6)\partial_1+(x_4x_5x_6+x_1x_3x_5)\partial_2\\
&+(x_1x_3x_4+x_4^2x_6)\partial_3+(x_1-x_3x_4x_6-x_1x_3^2)\partial_4\\
&+(-x_2x_4x_6-x_1x_2x_3)\partial_5+(-x_1x_4x_6-x_1^2x_3)\partial_6,
\end{split}
\end{equation*}
\begin{equation*}
\begin{split}
Y_{1,
5}&=(-x_5+x_5x_6^2+x_1x_2x_6)\partial_1+(x_5^2x_6+x_1x_2x_5)\partial_2\\
&+(x_4x_5x_6+x_1x_2x_4)\partial_3+(-x_3x_5x_6-x_1x_2x_3)\partial_4\\
&+(x_1-x_2x_5x_6-x_1x_2^2)\partial_5+(-x_1x_5x_6-x_1^2x_2)\partial_6,
\end{split}
\end{equation*}
\begin{equation*}
\begin{split}
Y_{1,
6}&=(-x_6+x_1^2x_6+x_6^3)\partial_1+(x_1^2x_5+x_6^2x_5)\partial_2\\
&+(x_1^2x_4+x_6^2x_4)\partial_3+(-x_1^2x_3-x_6^2x_3)\partial_4\\
&+(-x_1^2x_2-x_6^2x_2)\partial_5+(x_1-x_1^3-x_1x_6^2)\partial_6,
\end{split}
\end{equation*}
Thus, from \eqref{310} with $r=0$, \eqref{429} and \eqref{431},
taking account of Fig.\eqref{21}, we have
\begin{equation}\label{432}
\begin{split}
\phi 'Y_{1,
3}&=(x_1x_3x_6-x_1^2x_4)\partial_1+(x_2x_3x_6-x_1 x_2 x_4)\partial_2\\
&+(x_3^2x_6-x_1 x_3 x_4)\partial_3+(x_1+x_3x_4x_6-x_1x_4^2)\partial_4\\
&+(x_3x_5x_6-x_1x_4x_5)\partial_5+(-x_3+x_3x_6^2-x_1x_4x_6)\partial_6,
\end{split}
\end{equation}
\begin{equation*}
\begin{split}
\phi 'Y_{1,
4}&=(x_1x_4x_6+x_1^2x_3)\partial_1+(x_2x_4x_6+x_1x_2x_3)\partial_2\\
&+(-x_1+x_3x_4x_6+x_1x_3^2)\partial_3+(x_4^2x_6+x_1x_3x_4)\partial_4\\
&+(x_4x_5x_6+x_1x_3x_5)\partial_5+(-x_4+x_4x_6^2+x_1x_3x_6)\partial_6,
\end{split}
\end{equation*}
\begin{equation*}
\begin{split}
\phi 'Y_{1,
5}&=(x_1x_5x_6+x_1^2x_2)\partial_1+(-x_1+x_2x_5x_6+x_1x_2^2)\partial_2\\
&+(x_3x_5x_6+x_1x_2x_3)\partial_3+(x_4x_5x_6+x_1x_2x_4)\partial_4\\
&+(x_5^2x_6+x_1x_2x_5)\partial_5+(-x_5+x_5x_6^2+x_1x_2x_6)\partial_6,
\end{split}
\end{equation*}
\begin{equation*}
\begin{split}
\phi 'Y_{1,
6}&=(-x_1+x_1^3+x_6^2x_1)\partial_1+(x_1^2x_2+x_6^2x_2)\partial_2\\
&+(x_1^2x_3+x_6^2x_3)\partial_3+(x_1^2x_4+x_6^2x_4)\partial_4\\
&+(x_1^2x_5+x_6^2x_5)\partial_5+(-x_6+x_1^2x_6+x_6^3)\partial_6.
\end{split}
\end{equation*}
Thus, from \eqref{42}, \eqref{45} and \eqref{431}, we have
\begin{equation}\label{433}
\begin{split}
&d\eta(Y_{1, 2}, Y_{1, 3})=0,\quad d\eta(Y_{1, 2}, Y_{1, 4})=0,\\
&d\eta(Y_{1, 2}, Y_{1, 5})=-x_1^2,\quad d\eta(Y_{1, 2}, Y_{1,
6})=x_1x_2,\quad d\eta(Y_{1, b}, \xi)=0
\end{split}
\end{equation}
for any $b$ $(1 < b \leq 6)$. Similarly, from \eqref{422},
\eqref{431} and \eqref{432}, we have
\begin{equation}\label{434}
\begin{split}
&g(Y_{1, 2}, \phi 'Y_{1, 3})=0, \quad g(Y_{1, 2}, \phi 'Y_{1,
4})=0,\\
&g(Y_{1, 2}, \phi 'Y_{1, 5})=-x_1^2, \quad g(Y_{1, 2}, \phi 'Y_{1,
6})=x_1x_2, \quad g(Y_{1, b}, \phi '\xi)=0
\end{split}
\end{equation}
for any $b$ $(1 < b \leq 6)$. Therefore, from \eqref{13},
\eqref{433} and \eqref{434}, we can see that $(M_{r}, \phi ', \xi, \eta,
g)$ is a contact metric manifold. This also supports Theorem C. On the other hand, from (4.11) with $r=0$, $(M_{0}, \phi ', \xi,
\eta,g)$ is a space of constant sectional curvature 1. Therefore,
from the fact (\cite{B1}, Theorem 7.3), we see finally that
$(M_{0}, \phi ', \xi, \eta,g)$ is a Sasakian manifold. Further, taking account
of (3.1), we may check that, for each $r(-1<r<1)$, the map $F_{r}$ : $M_{r} \rightarrow M_{0}$
defined by $F_{r}(\sum_{i=1}^{6} x_{i}e_{i} + re_{7}) = \frac{1}{\sqrt{1-r^{2}}}(\sum_{i=1}^{6} x_{i}e_{i})$,
$(\sum_{i=1}^6 {x_{i}}^{2}=1-r^{2})$ on $M_{0}$ is a diffeomorphism from $M_{r}$ to $M_{0}$.
Thus, the pullback of the Sasakian structure $({\phi'}_{0}, {\xi}_{0},{\eta}_{0},{g}_{0})$ on $M_{0}$ to $M_{r}$
by the diffeomorphism $F_{r}$ is also a Sasakian structure on $M_{r}$ of constant sectional curvature 1 for each $r(-1<r<1)$.
We here note that the pull back Sasakian structure is given by $\bar{\phi}=({{F_{r}}^{-1}})_{*} \circ {\phi_{0}} \circ ({F_{r}})_{*}$,
$\bar{\xi} = ({{F_{r}}^{-1}})_{*}{\xi}_{0}=\xi$, $\bar{\eta}= {F_{r}}^{*}({\eta}_{0})=\eta$, $\bar{g}={F_{r}}^{*} (g_{0})=g$.
On the other hand, by modifying the above arguments suitably,
we may also check that $(M_{r}, \phi',\xi,\eta,g)$ is
not a contact metric manifold for any $r$ with $(-1<r<1, r\neq 0).$\\

\noindent {\bf{Remark 1.}} Let $M$ be a hypersurface in the nearly
K\"{a}hler unit 6-sphere $(S^6,J,\bar{g})$ oriented by unit normal
vector field $\nu$ and $(\phi,\xi,\eta,g)$ be the corresponding
naturally induced almost contact metric structure on $M$. Now, let
$G$ be the (1,2)-tensor field on $(S^6,J,\bar{g})$ given by
$G(\bar{X},\bar{Y}) = (\bar{\nabla}_{\bar{X}} J)(\bar{Y})$
for any $\bar{X}$,$\bar{Y}\in \mathfrak{X}(S^6)$, and $\psi$ be the
(1,1)-tensor field on $M$ defined by $\psi X = G(X,\nu)$ for any $X
\in\mathfrak{X}(M)$ \cite{DA}. Here, specifying $(M,\nu)$ as the
hypersurface $(M_{0},\nu)$ introduced in \S3, we can show that $\psi
= \phi'$ holds for $M_{0}$ by making use of the discussions in {\cite{DA, Se}}.\\\\
\noindent{\bf{Remark 2.}} J. Bernd, J. Bolton and L. M. Woodward have proved
that a Hopf hypersurface in the nearly {K\"ahler} 6-sphere $S^6$ is
either an open part of (i) a geodesic hypersphere of $S^6$ or (ii) a tube around an almost complex curve in $S^6$({\cite{BBW}},
Theorem 2). Taking account of this result, it seems also
meaningful to discuss the Hopf hypersurfaces of type (ii) in $S^6$
from the geometry of almost contact metric structures
viewpoint.

\section*{Acknowledgements}
\noindent This research was supported by Basic Science Research Program through
the National Research Foundation of Korea(NRF)
funded by the Ministry of Education (NRF-2016R1D1A1B03930449).


\medskip

\noindent

\end{document}